\magnification=\magstep1

\def\sg{\mathop{\rm sg}}
\def\Pf{\mathop{\rm Pf}}
\def\det{\mathop{\rm det}}

\def \sqr#1#2{{\vcenter{\hrule height.#2pt
    \hbox{\vrule width.#2pt height#1pt \kern#1pt \vrule width.#2pt}
    \hrule height.#2pt}}}

\centerline{An Application of Okada's Minor Summation Formula}
\centerline{to the Evaluation of a Multiple Integral}
\centerline{David P. Robbins}
\bigskip

\bigskip
\centerline{1.  INTRODUCTION}
\bigskip

Noam Elkies and Everett Howe [1] independently noticed 
a certain elegant product formula for the
multiple integral 
$$
\int_R\prod_{1\le i<j\le k} (x_j-x_i) dx_1\cdots dx_k,\eqno(1)
$$
where the region $R$ is the set of $k$-tuples $(x_1,\dots, x_k)$ satisfying
$ 0<x_1<\cdots <x_k < 1 .$   Later on Howe discovered that the formula was 
a special case of a formula for Selberg's integral described for
example in [2], page 339.

Here we prove an apparently different generalization
$$ 
\int_R \det\left(x_i^{a_j-1}\right) 
dx_1 \cdots dx_k = {\prod_{1 \le i<j \le k}(a_j-a_i)\over
\prod_{1 \le i \le k} a_i \prod_{1 \le i<j \le k} (a_j+a_i)}.\eqno(2)
$$
The original observation is the case $a_j=j$.  

The most interesting aspect of our proof is that we apply a limiting
form of a remarkable identity of Okada
[3] for summing the $k$ by $k$ minors of an $n$ by $k$ matrix.

First we describe Okada's formula.

Suppose that
$C$ is an $n$ by $k$ matrix.  For $1 \le i,j \le k$
define 
$$  S_{ij} = \sum_{1\le t<u \le n}(C_{ti}C_{uj}-C_{tj}C_{ui}).   $$
Note that $S$ is a skew-symmetric matrix.  $S_{ij}$ is the sum of 
all two-by-two minors of the submatrix of $C$ formed from columns $i$ and
$j$ (in that order).  

If $k$ is even, then we can form the Pfaffian of $S$, $\Pf(S)$, given
as
$$  
\sum \sg(\pi) S_{\pi(1),\pi(2)}S_{\pi(3),\pi(4)}\cdots
S_{\pi(k-1),\pi(k)}, 
$$
where the sum is over all permutations of $\{1,\dots,k\}$ satisfying
$$\pi(1)<\pi(2),\quad\pi(3)<\pi(4),\quad\ldots,\quad\pi(k-1)<\pi(k)$$
and 
$$ \pi(1) < \pi(3) < \cdots < \pi(k-1).  $$
It is known that the square of the Pfaffian is the
determinant of $S$.

Okada's formula states that the sum of the $k$ by $k$ minors of $C$ is
$\Pf(S)$.  It is easily proved.  Indeed it suffices to
prove the formula in the easy case of matrices $C$ with a single 1 in
every column.  Other cases follow since both sides of Okada's formula
are multilinear in the columns of $C$.

If $k$ is odd, we can form $S$ but not its Pfaffian.  However, Okada
observes that we can still sum the $k$ by
$k$ minors of $C$.  We simply adjoin a zeroth row and zeroth column to
$C$ which are both all zero except for a single 1 at their
intersection.  The sum of the $k+1$ by $k+1$ minors of this augmented
matrix is the same as the sum of the $k$ by $k$ minors of $C$.  So we
can apply the preceding case to the augmented matrix.  The effect on
$S$ is to add a zeroth row and a zeroth column with
$S_{0j}=-S_{j0}=\sum_{1 \le t \le n}C_{tj}$, the $j$th column sum of
$C$.  We can then sum the minors of $C$ by computing the Pfaffian of
the augmented $S$.  

Now let us return to the multiple integral.  
Notice that, for $ x_1 <\cdots < x_k$,  $\det \left(x_i^{a_j-1}\right)$ 
can be regarded as a minor of the $\infty$ by $k$ matrix
whose rows are indexed by $x$ with $0<x<1$ with the $x$th row equal to
$$
\matrix{
x^{a_1-1} &  x^{a_2-1} & \cdots &  x^{a_k-1}.
}
$$
Thus the integral can be regarded as a sort of sum of minors to which
Okada's formula can be applied.  

It is not difficult to make this idea more precise.  In fact there is
an integral form of Okada's identity.

\bigskip
\noindent{THEOREM.}
{\it 
Let $f_1,\dots,f_k$ be continuous functions on the interval $[a,b]$.
Let
$$ I= \int_{a<x_1<\cdots <x_k<b} \det(f_i(x_j))dx_1\cdots dx_k, $$
$$  I_{ij}=\int_{a<x<y<b}( f_i(x)f_j(y)-f_j(x)f_i(y))dx dy, $$
and
$$  I_i=\int_{a<x<b}f_i(x) dx.  $$
If $k$ is even, then $I$ is the Pfaffian of the $k$ by $k$ matrix 
$I_{ij}$, $1\le i,j \le k$.  If $k$ is odd, then $I$ is the
Pfaffian of the $k+1$ by $k+1$ matrix which is $I_{ij}$ augmented with
a zeroth row and column defined by $I_{00}=0$ and 
$I_{0i}=-I_{i0}=I_i$ for $1\le i \le k$.
}

\noindent{\it Proof}:
Let $n$ be an integer, $h=(b-a)/n$ and $x_i=a+hi$, 
$ i=0,\dots,n$.  Consider the matrix  
$$C=h f_j(x_i)_{0 \le i \le n, 1 \le j \le k}.$$
Let $S$ be the sum of the $k$ by $k$ minors of $C$,
let $S_{ij}$ be the sum of the two by two minors of
the matrix formed from columns $i$ and $j$ of $C$
(in that order) and let $S_i$ be the sum of the entries
of column $i$ of $C$.  $S$, $S_i$, and $S_{ij}$ all
depend on $n$.
Standard properties of integrals yield
$$  I=\lim_{n \rightarrow \infty} S, \quad  I_i=\lim_{n \rightarrow \infty}S_i, \quad 
 I_{ij}=\lim_{n \rightarrow \infty}S_{ij}. $$
The rest then follows from Okada's formula.

\medskip

We can now apply Theorem~1
to prove (2).
Let
$$I=\int_R \det\left(x_i^{a_j-1}\right) dx_1 \cdots dx_k. $$
Then routine integration shows that
$$
I_i = {1\over a_i}, \quad  I_{ij}={(a_j-a_i)\over a_ia_j(a_i+a_j)}. \eqno{(3)}
$$

From Theorem~1, if $k$ is even, $I$ is the 
Pfaffian of the $k$ by $k$ matrix 
$$
\left[\matrix{
0 & {a_2-a_1\over a_1a_2(a_1+a_2)} & {a_3-a_1\over a_1a_3(a_1+a_3)}& 
\cdots & {a_k-a_1\over a_1a_k(a_1+a_k)}\cr
{a_1-a_2\over a_1a_2(a_1+a_2)} & 0 & {a_3-a_2\over a_2a_3(a_2+a_3)}&
\cdots\cr
{a_1-a_3\over a_1a_3(a_1+a_3)} & {a_2-a_3\over a_2a_3(a_2+a_3)} &
0\cr
\vdots & \vdots & &\ddots\cr
 {a_1-a_k\over a_1a_k(a_1+a_k)} & & & & 0\cr
}\right]
$$
while if $k$ is odd, $I$ is the Pfaffian of the $k+1$ by $k+1$ matrix
$$
\left[\matrix{
0 & {1\over a_1} & {1\over a_2} & \cdots & & {1\over a_k}\cr
-{1\over a_1} & 0 & {a_2-a_1\over a_1a_2(a_1+a_2)} &\cdots & &
{a_k-a_1\over a_1a_k(a_1+a_k)}\cr
 -{1\over a_2} & {a_1-a_2\over a_1a_2(a_1+a_2)} & 0 \cr
\vdots & \vdots & &\ddots\cr
 {a_1\over a_k} & {a_1-a_k\over a_1a_k(a_1+a_k)} & & & &0
}\right].
$$
In each case the determinant is a homogeneous rational function of
degree $-2k$ so the Pfaffian is a homogeneous rational function of
degree $-k$.  However, since $I$ vanishes when any two of the $a$'s
are equal, we know that the numerator of $I$ must contain a factor of
$N=\prod_{1 \le i<j \le k}(a_j-a_i) $.  Hence the denominator is of
degree $k+k(k-1)/2+e$, where $e$ is the degree of the remaining factor
of the numerator.  On the other hand one also sees directly from the
Pfaffian expressions that the denominator must divide $D=\prod_{1 \le
i \le k} a_i \prod_{1\le i<j\le k} (a_i+a_j).$ It follows that $e$ is
zero and that the numerator of $I$ is $N$ to within a constant factor
and the denominator $D$ to within a constant factor.  Thus the left
side of (2) is equal to the right to within a constant factor.

We now need to show that this constant factor is 1.  This is certainly
true for $k=1$ and $k=2$ by (3).  For larger $k$, when $k$
is even, if we multiply the first row and column of the Pfaffian
matrix by $a_1$ and then set $a_1=0$ we obtain the matrix for the
case $k-1$ (applied to $a_2,\dots,a_k$).  Similarly when $k$ is odd,
if we multiply the second row and second column of the Pfaffian by
$a_1$, we can easily reduce to the case $k-1$.  On the other hand the
right side of (2) also has the property that when we multiply
it by $a_1$ and then set $a_1=0$ we obtain the right side of
(2) for $k-1$ (applied to $a_2,\dots,a_k$).  These
observations yield an immediate inductive proof that the constant is
1.

\bigskip
\centerline{REFERENCES}
\bigskip

\item{1.}
Everett Howe, Oral communication.

\item{2.}
Madan Lal Mehta, Random Matrices, {\it Academic Press}, 1991.  

\item{3.} 
Soichi Okada, {\it On the Generating Functions for Certain Classes of
Plane Partitions}, Journal of Combinatorial Theory, Volume 1, Series
A, pages 1--23, 1989.

\bye